\documentclass[a4paper,12pt,twoside]{article}
\usepackage{amsmath,amssymb,ifthen}
\usepackage[amsmath,thmmarks]{ntheorem}
\usepackage{graphicx}
\usepackage{paper}
\usepackage[all]{xy}
\usepackage{mathrsfs}
\usepackage{url}

\DeclareMathOperator{\Sw}{Sw}
\DeclareMathOperator{\Ar}{Ar}
\DeclareMathOperator{\Kl}{Kl}
\DeclareMathOperator{\tr}{tr}
\DeclareMathOperator{\mult}{mult}
\DeclareMathOperator{\add}{add}
\DeclareMathOperator{\Sp}{Sp}
\DeclareMathOperator{\SO}{SO}
\DeclareMathOperator{\SL}{SL}
\DeclareMathOperator{\St}{\mathbf{St}}
\DeclareMathOperator{\Ad}{Ad}
\DeclareMathOperator{\Cent}{Cent}
\DeclareMathOperator{\Lie}{Lie}
\DeclareMathOperator{\cInd}{c-Ind}
\DeclareMathOperator{\Sym}{Sym}
\newcommand{\G}{\mathbb{G}}

\SelectTips{cm}{11}
\title{On the formal degree conjecture for simple supercuspidal representations}
\author{Yoichi Mieda}

\begin{document}
\maketitle

\begin{firstfootnote}
 Graduate School of Mathematical Sciences, The University of Tokyo, 3--8--1 Komaba, Meguro-ku, Tokyo, 153--8914, Japan

 E-mail address: \texttt{mieda@ms.u-tokyo.ac.jp}

\end{firstfootnote}

\begin{abstract}
 We prove the formal degree conjecture for simple supercuspidal representations of symplectic groups and quasi-split
 even special orthogonal groups over a $p$-adic field,
 under the assumption that $p$ is odd. The essential part is to compute the Swan conductor of the exterior square of
 an irreducible local Galois representation with Swan conductor $1$.
 It is carried out by passing to the equal characteristic local field and using the theory of Kloosterman sheaves.
\end{abstract}

\section{Introduction}
Let $F$ be a $p$-adic field, and $G$ a connected reductive group over $F$.
For an irreducible discrete series representation $\pi$ of $G(F)$, we can consider an invariant $\deg(\pi)\in\R_{>0}$
called the formal degree of $\pi$. It is in some sense a generalization of the dimension (or the degree)
of a finite-dimensional representation.
On the other hand, by the local Langlands correspondence,
irreducible smooth representations of $G(F)$ are conjecturally parametrized by pairs $(\phi,\rho)$,
where $\phi\colon W_F\times\SL_2(\C)\to {}^LG$ is an $L$-parameter,
and $\rho$ is an irreducible representation of a finite group $\mathcal{S}_\phi$ determined by $\phi$.
The formal degree conjecture, which was proposed by Hiraga-Ichino-Ikeda \cite{MR2350057},
predicts that $\deg(\pi)$ can be described
by using the pair $(\phi_\pi,\rho_\pi)$ attached to $\pi$. For more precise formulation, see Section \ref{sec:FDC-ssc}.
This conjecture has been solved for general linear groups \cite{MR2350057},
odd special orthogonal groups \cite{MR3649356} and unitary groups \cite{2018arXiv181200047B},
but it seems still open for many other groups.

In this article, we will focus on a very special class of discrete series representations, 
simple supercuspidal representations.
They are introduced in \cite{MR2730575} and \cite{MR3164986}, and characterized among irreducible smooth representations
by the property that they have minimal positive depth.
The local Langlands correspondence for simple supercuspidal representations of quasi-split classical groups has been investigated in
Oi's work \cite{Oi-ssc-classical} very precisely. As an application of his results, he obtained the following theorem:

\begin{thm}[{{\cite[Theorem 9.3]{Oi-ssc-classical}}}]\label{thm:FDC-Oi-intro}
 Let $n\ge 1$ be an integer and write $2n=p^en'$ with $p\nmid n'$.
 Assume $p\neq 2$ and either $p\nmid 2n$ or $n'\mid p-1$.
 Let $G$ be one of the following groups:
 \begin{itemize}
  \item $\Sp_{2n}$,
  \item the quasi-split $\SO_{2n}$ attached to a ramified quadratic extension of $F$,
  \item the split $\SO_{2n+2}$,
  \item or the quasi-split $\SO_{2n+2}$ attached to an unramified quadratic extension of $F$.
 \end{itemize}
 Then, the formal degree conjecture holds for simple supercuspidal representations of $G(F)$.
\end{thm}

The goal of this article is to remove the condition ``either $p\nmid 2n$ or $n'\mid p-1$'' in the theorem above.
Here is our main theorem:

\begin{thm}[Theorem \ref{thm:main} and Remark \ref{rem:Gan-Ichino}]\label{thm:main-intro}
 Assume $p\neq 2$. Let $G$ be one of the groups in Theorem \ref{thm:FDC-Oi-intro}.
 Then, the formal degree conjecture holds for simple supercuspidal representations of $G(F)$.
\end{thm}

By the same method as in \cite{Oi-ssc-classical}, Theorem \ref{thm:main-intro} is easily reduced to the following:

\begin{thm}[Theorem \ref{thm:exterior-square}]\label{thm:exterior-square-intro}
 Let $\tau$ be a $2n$-dimensional irreducible smooth representation of $W_F$ such that $\Sw \tau=1$.
 Then we have $\Sw(\wedge^2\tau)=n-1$.
\end{thm}

In the case $p\mid 2n$ and $n'\mid p-1$ (which is more difficult than the case $p\nmid 2n$),
Oi used an explicit description of $\tau$ in \cite{2015arXiv150902960I} to obtain the theorem above.
Our strategy to Theorem \ref{thm:exterior-square-intro} is totally different.
First we use Deligne's result \cite{MR771673} to reduce Theorem \ref{thm:exterior-square-intro} to the case
where $F$ is an equal characteristic local field.
In the equal characteristic case, every irreducible smooth representation of $W_F$ with Swan conductor $1$ is essentially
obtained as the localization at $\infty\in \P^1$ of a Kloosterman sheaf $\Kl$
(see \cite[Sommes.\ trig.]{MR0463174} and \cite{MR955052}). The Swan conductor of the localization at $\infty$
of $\wedge^2\Kl$ can be computed by using
the Grothendieck-Ogg-Shafarevich formula and the Grothendieck-Lefschetz trace formula.

The outline of this paper is as follows. 
In Section \ref{sec:exterior}, we will compute the Swan conductor of the exterior square of
an irreducible smooth representation $\tau$ of $W_F$ with Swan conductor $1$.
Although in Theorem \ref{thm:exterior-square-intro} we assumed that $\dim \tau$ is even
(in fact only this case is needed to prove Theorem \ref{thm:main-intro}),
we will also treat the case where $\dim \tau$ is odd.
In Section \ref{sec:FDC-ssc}, after recalling the formal degree conjecture,
we deduce Theorem \ref{thm:main-intro} from Theorem \ref{thm:exterior-square-intro}.

\medbreak
\noindent{\bfseries Acknowledgment}\quad
This work was supported by JSPS KAKENHI Grant Number 15H03605.

\medbreak
\noindent{\bfseries Notation}\quad
Every representation is considered over $\C$, unless otherwise noted. 

\section{Exterior square of local Galois representations with Swan conductor 1}\label{sec:exterior}
Let $p$ be a prime number and $F$ a finite extension of $\Q_p$. We write $k$ for the residue field of $F$
and $q$ for the cardinality of $k$.
We fix an algebraic closure $\overline{F}$ of $F$
and put $\Gamma_F=\Gal(\overline{F}/F)$. Let $W_F$ denote the Weil group of $F$, that is,
the subgroup of $\Gamma_F$ consisting of elements which induce integer powers of the Frobenius automorphism on 
the residue field $\overline{k}$ of $\overline{F}$.
It is a locally compact group containing the inertia group $I_F$ as an open subgroup.

Recall that $\Gamma_F$ is equipped with the upper numbering ramification filtration $\{\Gamma_F^j\}_{j\in \R_{\ge 0}}$,
which is a descending filtration consisting of open normal subgroups of $\Gamma_F$.
The subgroup $\Gamma_F^0$ equals the inertia group $I_F$, and $\Gamma_F^{0+}$ equals the wild inertia group $P_F$.
Here $\Gamma_F^{j+}$ denotes the closure of $\bigcup_{j'>j}\Gamma_F^{j'}$, as usual.

Let $V$ be a finite-dimensional smooth representation of $W_F$.
It is known that there exists a unique direct sum decomposition
$V=\bigoplus_{j\in\R_{\ge 0}}V_j$ as a representation of $P_F$, called the break decomposition,
such that
\begin{itemize}
 \item $V_0=V^{P_F}$, and
 \item $V_j^{\Gamma_F^j}=0$ and $V_j^{\Gamma_F^{j+}}=V_j$ for each $j\in\R_{>0}$
\end{itemize}
(see \cite[Proposition 1.1, Lemma 1.4]{MR955052}).
The numbers $j$ with $V_j\neq 0$ are called the breaks of $V$. The Swan conductor $\Sw V$ of $V$ is defined by
\[
 \Sw V=\sum_{j\in \R_{\ge 0}}j\dim V_j.
\]
Note that $\Sw V$ depends only on the restriction of $V$ to $P_F$.

Let $n\ge 1$ be an integer. In this section, we prove the following result.

\begin{thm}\label{thm:exterior-square}
 Let $(\tau,V)$ be an $n$-dimensional irreducible smooth representation of $W_F$ such that $\Sw \tau=1$.
 Then we have
 \[
  \Sw(\wedge^2\tau)=\begin{cases}
		     m-1 & \text{if $n=2m$ is even,}\\
		     m & \text{if $n=2m+1$ is odd.}
		    \end{cases}
 \]
\end{thm}

If $n=1$, this theorem is obvious. Therefore, we assume $n\ge 2$ in the following.
First we notice a simple lemma.

\begin{lem}\label{lem:totally-wild}
 Let $(\tau,V)$ be as in Theorem \ref{thm:exterior-square}.
 Then we have $V^{P_F}=0$. Moreover, $V$ has only one break $1/n$ and $V\vert_{I_F}$ is irreducible.
\end{lem}

\begin{prf}
 Since $P_F$ is a normal subgroup of $W_F$, $V^{P_F}$ is a $W_F$-subrepresentation of $V$.
 The condition $\Sw V=1$ implies that $V^{P_F}\neq V$. Therefore we have $V^{P_F}=0$ by the irreducibility of $V$.
 By \cite[Lemma 1.11]{MR955052}, $V$ has only one break $1/n$ and $V\vert_{I_F}$ is irreducible.
\end{prf}

To prove Theorem \ref{thm:exterior-square}, we will pass to the equal characteristic case.
Put $F'=k((T))$, which is an equal characteristic local field.
By Deligne's result \cite{MR771673}, we can prove the following:

\begin{lem}\label{lem:passage-to-char-p}
 Let $(\tau,V)$ be as in Theorem \ref{thm:exterior-square}.
 Then there exists an $n$-dimensional irreducible smooth representation $(\tau',V')$ of $W_{F'}$ such that
 \[
  \Sw \tau'=\Sw \tau=1,\qquad \Sw(\wedge^2\tau')=\Sw(\wedge^2\tau).
 \]
\end{lem}

\begin{prf}
 By \cite[\S 3.5]{MR771673}, there exists an isomorphism $\Gamma_F/\Gamma_F^1\cong \Gamma_{F'}/\Gamma_{F'}^1$, which is canonical up to
 inner automorphisms. By construction, it preserves the upper numbering ramification filtrations of $\Gamma_F$ and $\Gamma_{F'}$.
 Further, \cite[Proposition 3.6.1]{MR771673} tells us that it induces an isomorphism $W_F/\Gamma_F^1\cong W_{F'}/\Gamma_{F'}^1$.
 By using this isomorphism, we can construct a functor
 \begin{itemize}
  \item from the category of finite-dimensional smooth representations of $W_F$
	whose breaks are less than $1$
  \item to the category of finite-dimensional smooth representations of $W_{F'}$
	whose breaks are less than $1$.
 \end{itemize}
 Clearly this functor maps irreducible representations to irreducible representations,
 commutes with exterior products, and preserves the Swan conductors.

 By Lemma \ref{lem:totally-wild}, $\tau$ has only one break $1/n$, which is less than $1$.
 Therefore we can take $(\tau',V')$ as the image of $(\tau,V)$ under this functor.
\end{prf}

Now we use the Kloosterman sheaves introduced in \cite[Sommes.\ trig.]{MR0463174} and \cite{MR955052}.
Let us recall their construction briefly.
Take a prime number $\ell\neq p$. We fix an isomorphism $\overline{\Q}_\ell \cong\C$ and identify them.
Note that an irreducible finite-dimensional continuous representation of $W_F$ over $\overline{\Q}_\ell$ is automatically smooth,
hence can be identified with an irreducible smooth representation of $W_F$ over $\C$.
Let $\mathbb{P}^1$ denote the projective line over $k$, and put $\A^1=\P^1\setminus \{\infty\}$, $\G_m=\P^1\setminus \{0,\infty\}$.
We consider the diagram
\[
 \G_m\xleftarrow{\mult}\G_m^n\xrightarrow{\add}\A^1,
\]
where the maps $\mult$ and $\add$ are given by $(x_1,\ldots,x_n)\mapsto x_1\cdots x_n$ and 
$(x_1,\ldots,x_n)\mapsto x_1+\cdots+x_n$, respectively.
We fix a non-trivial additive character $\psi\colon k\to \C^\times$, and
write $\mathcal{L}_\psi$ for the Artin-Schreier sheaf on $\A^1$ corresponding to $\psi$.
For multiplicative characters $\chi_1,\ldots,\chi_n\colon k^\times\to \C^\times$, we can construct the Kummer sheaf
$\mathcal{K}_{\chi_1},\ldots,\mathcal{K}_{\chi_n}$ on $\G_m$. We put
\[
 \Kl(\chi_1,\ldots,\chi_n)=R\mult_!\bigl((\mathcal{K}_{\chi_1}\boxtimes\cdots\boxtimes \mathcal{K}_{\chi_n})\otimes \add^*\mathcal{L}_\psi\bigr)[-n+1].
\]
If $\chi_1=\cdots=\chi_n=1$, we simply write $\Kl_n$ for $\Kl(\chi_1,\ldots,\chi_n)$.
It is known that $\Kl(\chi_1,\ldots,\chi_n)$ is a smooth sheaf on $\G_m$ of rank $n$.
Further, it enjoys the following properties:
\begin{itemize}
 \item $\Kl(\chi_1,\ldots,\chi_n)_{\overline{0}}$, which is a representation of $\Gamma_{\Frac \widehat{\mathcal{O}}_{\P^1,0}}$,
       is tamely ramified.
 \item $\Kl(\chi_1,\ldots,\chi_n)_{\overline{\infty}}$,
       which is a representation of $\Gamma_{\Frac \widehat{\mathcal{O}}_{\P^1,\infty}}$,
       is totally wildly ramified with Swan conductor $1$ (in particular it is irreducible by \cite[Lemma 1.11]{MR955052}).
\end{itemize}
Here $\widehat{\mathcal{O}}_{\P^1,x}$ denotes the completion of the local ring $\mathcal{O}_{\P^1,x}$ at $x\in \P^1$.
See \cite[Sommes.\ trig., Th\'eor\`eme 7.8]{MR0463174} and \cite[Theorem 4.1.1]{MR955052} for detail.

In the following, we fix an isomorphism $k[[T]]\cong \widehat{\mathcal{O}}_{\P^1,\infty}$ and identify them.
Then $\Kl(\chi_1,\ldots,\chi_n)_{\overline{\infty}}$ can be regarded as
an $n$-dimensional irreducible smooth representation of $W_{F'}$.

\begin{lem}\label{lem:tau'-Kloosterman}
 Let $\tau'$ be an $n$-dimensional irreducible smooth representation of $W_{F'}$ with Swan conductor $1$.
 Then we have $\Sw(\wedge^2\tau')=\Sw(\wedge^2\Kl_{n,\overline{\infty}})$.
\end{lem}

\begin{prf}
 By replacing $\tau'$ by its unramified twist,
 we may assume that $\tau'$ extends to a smooth representation of $\Gamma_{F'}$.
 Note that $\tau'$ is defined over a finite extension $E_\lambda$ of $\Q_\ell$ contained in $\overline{\Q}_\ell$.
 By the theorem of Katz-Gabber (\cite[Theorem 1.5.6]{MR867916}),
 there exists a smooth $E_\lambda$-sheaf $\mathcal{F}$ on $\G_m$ of rank $n$ such that
 $\mathcal{F}_{\overline{0}}$ is tamely ramified and $\mathcal{F}_{\overline{\infty}}$ is isomorphic to $\tau'$.
 By \cite[Theorem 8.7.1]{MR955052} (see also the proof of \cite[Corollary 8.7.2]{MR955052}),
 there exist a finite extension $k'$ of $k$, an element $a'\in k'^\times$ and multiplicative characters
 $\chi'_1,\ldots,\chi'_n\colon k'^\times\to \C^\times$ such that
 \[
  \mathcal{F}\otimes_kk'\cong \iota_{a'}^*\Kl(\chi'_1,\ldots,\chi'_n),
 \]
 where $\iota_{a'}\colon \G_m\otimes_kk'\to \G_m\otimes_kk'$ is the multiplication by $a'$ and $\Kl(\chi'_1,\ldots,\chi'_n)$
is the Kloosterman sheaf over $\G_m\otimes_kk'$ with respect to the additive character $\psi\circ\tr_{k'/k}$ of $k'$.
 Since the base change from $k$ to $k'$ and the pull-back by $\iota_{a'}$ do not affect the Swan conductor at $\infty$,
 we conclude that $\Sw(\wedge^2\tau')=\Sw(\wedge^2\Kl(\chi'_1,\ldots,\chi'_n)_{\overline{\infty}})$.
 On the other hand, by \cite[Proposition 10.1]{MR955052},
 the restriction of $\Kl(\chi'_1,\ldots,\chi'_n)_{\overline{\infty}}$ to $P_{k'((T))}=P_{F'}$ is
 independent of $\chi'_1,\ldots,\chi'_n$. 
 Hence we have 
 \begin{align*}
  \Sw(\wedge^2\Kl(\chi'_1,\ldots,\chi'_n)_{\overline{\infty}})&=\Sw(\wedge^2\Kl(1',\ldots,1')_{\overline{\infty}})=\Sw(\wedge^2(\Kl_n\otimes_kk')_{\overline{\infty}})\\
  &=\Sw(\wedge^2\Kl_{n,\overline{\infty}}),
 \end{align*}
 where $1'$ denotes the trivial character of $k'^\times$. This concludes the proof.
\end{prf}

By Lemma \ref{lem:tau'-Kloosterman}, we may focus on computing $\Sw(\wedge^2\Kl_{n,\overline{\infty}})$.
Since $\wedge^2\Kl_{n,\overline{0}}$ is tame,
the Grothendieck-Ogg-Shafarevich formula \cite[Expos\'e X, Th\'eor\`eme 7.1]{SGA5} tells us that
\[
 \Sw(\wedge^2\Kl_{n,\overline{\infty}})=-\chi_c(\G_m,\wedge^2\Kl_n)
 :=-\sum_{i=0}^2 (-1)^i\dim H^i_c(\G_m\otimes_k\overline{k},\wedge^2\Kl_n).
\]
We shall determine the Euler characteristic $\chi_c(\G_m,\wedge^2\Kl_n)$ by computing the $L$-function 
\[
 L(\G_m,\wedge^2\Kl_n,X)=\exp\biggl(\sum_{r=1}^\infty \Bigl(\sum_{a\in k_r^\times}\Tr(\Frob_a,\wedge^2\Kl_{n,\overline{a}})\Bigr)\frac{X^r}{r}\biggr),
\]
where $k_r=\F_{q^r}$ denotes the degree $r$ extension of $k=\F_q$.

\begin{prop}\label{prop:exterior-L-function}
 We have
 \[
  L(\G_m,\wedge^2\Kl_n,X)=\begin{cases}
			   \dfrac{(1-qX)(1-q^3X)\cdots (1-q^{2m-1}X)}{1-q^{2m}X}& \text{if $n=2m$ is even,}\\
			   (1-qX)(1-q^3X)\cdots (1-q^{2m-1}X) & \text{if $n=2m+1$ is odd.}
			  \end{cases}
 \]
\end{prop}

\begin{prf}
 First note that
 \[
  \sum_{a\in k_r^\times}\Tr(\Frob_a,\wedge^2\Kl_{n,\overline{a}})
 =\frac{1}{2}\Bigl(\sum_{a\in k_r^\times}\Tr(\Frob_a,\Kl_{n,\overline{a}}\otimes \Kl_{n,\overline{a}})-\sum_{a\in k_r^\times}\Tr(\Frob_a^2,\Kl_{n,\overline{a}})\Bigr).
 \]

 By the proof of \cite[Proposition 10.4.1]{MR955052}, we have
 \begin{align*}
  &\sum_{a\in k_r^\times}\Tr(\Frob_a,\Kl_{n,\overline{a}}\otimes \Kl_{n,\overline{a}})=S_r(n,1,1)\\
  &\qquad=\begin{cases}
     -1-q^r-q^{2r}-\cdots-q^{(n-1)r}+q^{nr}& \text{if $n$ is even or $p=2$,}\\
     -1-q^r-q^{2r}-\cdots-q^{(n-1)r}& \text{if $n$ is odd and $p\neq 2$}
    \end{cases}
 \end{align*}
 (note that ``if $\alpha(-\alpha)^n\beta=1$'' in the end of p.~173 of \cite{MR955052} should be ``if $\alpha(-\alpha)^n\beta=-1$'').
 On the other hand, by \cite[(4.2.1.3), (4.2.1.5)]{MR955052}, we have
 \[
 \sum_{a\in k_r^\times}\Tr(\Frob_a^2,\Kl_{n,\overline{a}})=(-1)^{n-1}\frac{1}{q^{2r}-1}\sum_{\rho\in(k_{2r}^\times)^\vee} g(\psi\circ \tr_{k_{2r}/k},\rho^{q^r-1})^n,
 \]
 where $(k_{2r}^\times)^\vee$ denotes the set of characters of $k_{2r}^\times$ and 
 \[
  g(\psi\circ \tr_{k_{2r}/k},\rho^{q^r-1})=\sum_{a\in k_{2r}^\times}\psi(\tr_{k_{2r}/k}(a))\rho(a)^{q^r-1}
 \]
 denotes the Gauss sum.
 Further, by \cite[(4.2.1.13)]{MR955052}, we have
 \[
  g(\psi\circ \tr_{k_{2r}/k},\rho^{q^r-1})
  =\begin{cases}
    q^r\rho(-1)& \text{if $\rho^{q^r-1}\neq 1$,}\\
    -1 & \text{if $\rho^{q^r-1}=1$.}
   \end{cases}
 \]
 Since 
 \begin{align*}
  \#\{\rho\in(k_{2r}^\times)^\vee\mid \rho^{q^r-1}=1\}&=q^r-1,\\
  \#\{\rho\in(k_{2r}^\times)^\vee\mid \rho^{q^r-1}\neq 1,\rho(-1)=1\}
  &=\begin{cases}
    \dfrac{(q^r-1)^2}{2} &\text{if $p\neq 2$,}\\
    q^{2r}-q^r &\text{if $p=2$,}
   \end{cases}\\
  \#\{\rho\in(k_{2r}^\times)^\vee\mid \rho^{q^r-1}\neq 1,\rho(-1)=-1\}
   &=\begin{cases}
    \dfrac{q^{2r}-1}{2} &\text{if $p\neq 2$,}\\
    0 &\text{if $p=2$,}
    \end{cases}
 \end{align*}
 we have
 \[
  \sum_{a\in k_r^\times}\Tr(\Frob_a^2,\Kl_{n,\overline{a}})
 =\begin{cases}
   \dfrac{1}{q^r+1}\Bigl(\dfrac{(-1)^{n-1}q^{nr}(q^r-1)}{2}-\dfrac{q^{nr}(q^r+1)}{2}-1\Bigr) &\text{if $p\neq 2$,}\\[10pt]
   \dfrac{(-1)^{n-1}q^{(n+1)r}-1}{q^r+1} &\text{if $p=2$.}
   \end{cases}
 \]
 Now we assume that $n=2m$ is even. Then we have
 \[
  \sum_{a\in k_r^\times}\Tr(\Frob_a^2,\Kl_{n,\overline{a}})=-\frac{q^{(n+1)r}+1}{q^r+1}=-1+q^r-q^{2r}+\cdots-q^{2mr},
 \]
 hence
 \begin{align*}
  &\sum_{a\in k_r^\times}\Tr(\Frob_a,\wedge^2\Kl_{n,\overline{a}})\\
  &\qquad =\frac{1}{2}\Bigl((-1-q^r-q^{2r}-\cdots-q^{(2m-1)r}+q^{2mr})\\
  &\qquad\qquad\qquad\qquad\qquad\qquad\qquad -(-1+q^r-q^{2r}+\cdots+q^{(2m-1)r}-q^{2mr})\Bigr)\\
  &\qquad =-q^r-q^{3r}-\cdots-q^{(2m-1)r}+q^{2mr}.
 \end{align*}
 Therefore we conclude that
 \[
  L(\G_m,\wedge^2\Kl_n,X)=\frac{(1-qX)(1-q^3X)\cdots(1-q^{2m-1}X)}{1-q^{2m}X}.
 \]
 Next we consider the case where $n=2m+1$ is odd and $p\neq 2$. We have
 \[
  \sum_{a\in k_r^\times}\Tr(\Frob_a^2,\Kl_{n,\overline{a}})=-\frac{q^{nr}+1}{q^r+1}=-1+q^r-q^{2r}+\cdots-q^{2mr},
 \]
 hence
 \begin{align*}
  &\sum_{a\in k_r^\times}\Tr(\Frob_a,\wedge^2\Kl_{n,\overline{a}})\\
  &\qquad =\frac{1}{2}\Bigl((-1-q^r-q^{2r}-\cdots-q^{(2m-1)r}-q^{2mr})\\
  &\qquad\qquad\qquad\qquad\qquad\qquad\qquad -(-1+q^r-q^{2r}+\cdots+q^{(2m-1)r}-q^{2mr})\Bigr)\\
  &\qquad =-q^r-q^{3r}-\cdots-q^{(2m-1)r}.
 \end{align*}
 Therefore we conclude that
 \[
  L(\G_m,\wedge^2\Kl_n,X)=(1-qX)(1-q^3X)\cdots(1-q^{2m-1}X).
 \]
 Finally we assume that $n=2m+1$ is odd and $p=2$. 
 Then we have
 \[
  \sum_{a\in k_r^\times}\Tr(\Frob_a^2,\Kl_{n,\overline{a}})=\frac{q^{(n+1)r}-1}{q^r+1}=-1+q^r-q^{2r}+\cdots+q^{(2m+1)r},
 \]
 hence
 \begin{align*}
  &\sum_{a\in k_r^\times}\Tr(\Frob_a,\wedge^2\Kl_{n,\overline{a}})\\
  &\qquad =\frac{1}{2}\Bigl((-1-q^r-q^{2r}-\cdots-q^{(2m-1)r}-q^{2mr}+q^{(2m+1)r})\\
  &\qquad\qquad\qquad\qquad\qquad\qquad\qquad -(-1+q^r-q^{2r}+\cdots-q^{2mr}+q^{(2m+1)r})\Bigr)\\
  &\qquad =-q^r-q^{3r}-\cdots-q^{(2m-1)r}.
 \end{align*}
 Therefore we conclude that
 \[
  L(\G_m,\wedge^2\Kl_n,X)=(1-qX)(1-q^3X)\cdots(1-q^{2m-1}X).
 \]
\end{prf}

By the Grothendieck-Lefschetz trace formula, we have
\[
 L(\G_m,\wedge^2\Kl_n,X)=\prod_{i=0}^2 \det(1-X\Frob;H^i_c(\G_m\otimes_k\overline{k},\wedge^2\Kl_n))^{(-1)^{i+1}}.
\]
In particular we have $\deg L(\G_m,\wedge^2\Kl_n,X)=-\chi_c(\G_m,\wedge^2\Kl_n)$.
Hence we obtain the following corollary:

\begin{cor}\label{cor:exterior-Euler-char}
 We have
 \[
  \Sw(\wedge^2\Kl_{n,\overline{\infty}})=-\chi_c(\G_m,\wedge^2\Kl_n)=\begin{cases}
		     m-1 & \text{if $n=2m$ is even,}\\
		     m & \text{if $n=2m+1$ is odd.}
			  \end{cases}
 \]
\end{cor}

Now Theorem \ref{thm:exterior-square} follows from Lemmas \ref{lem:passage-to-char-p}, \ref{lem:tau'-Kloosterman}
and Corollary \ref{cor:exterior-Euler-char}.

\section{The formal degree conjecture for simple supercuspidal representations}\label{sec:FDC-ssc}
In this section, we deduce the formal degree conjecture for simple supercuspidal representations of symplectic groups
and quasi-split even special orthogonal groups from Theorem \ref{thm:exterior-square}.
Let us first recall the conjecture quickly in the case of symplectic groups.
For more detail, see \cite{MR2350057}.

In the following, we put $G=\Sp_{2n}$ for an integer $n\ge 2$.
We fix a non-trivial additive character $\psi\colon F\to \C^\times$
and a Haar measure on $G(F)$.
Let $(\pi,V)$ be an irreducible discrete series representation of $G(F)$.
We fix a $G(F)$-invariant inner product $(\ ,\ )\colon V\times V\to \C$.
Then, there exists a unique positive real number $\deg(\pi)$, the formal degree of $\pi$,
satisfying
\[
 \int_{G(F)} (\pi(g)v,w)\overline{(\pi(g)v',w')}dg=\deg(\pi)^{-1}(v,v')\overline{(w,w')}
\]
for every $v,w,v',w'\in V$. It depends on the fixed measure on $G(F)$, but is independent of the inner product $(\ ,\ )$.
The formal degree conjecture predicts that $\deg(\pi)$ can be described by using the local Langlands correspondence.

By the local Langlands correspondence due to Arthur \cite{MR3135650}, discrete series representations of $G(F)$
are parametrized by pairs $(\phi,\rho)$, where
\begin{itemize}
 \item $\phi\colon W_F\times\SL_2(\C)\to \widehat{G}(\C)=\SO_{2n+1}(\C)$ is an $L$-parameter
       such that the centralizer group $S_\phi=\Cent_{\widehat{G}(\C)}(\Imm \phi)$ is finite,
 \item and $\rho$ is an irreducible representation of $\pi_0(S_\phi)=S_\phi$.
\end{itemize}
The pair attached to a discrete series representation $\pi$ is denoted by $(\phi_\pi,\rho_\pi)$.

Here is the statement of the formal degree conjecture for $\Sp_{2n}$:

\begin{conj}[{{\cite[Conjecture 1.4]{MR2350057}}}]\label{conj:HII}
 For an irreducible discrete series representation $\pi$ of $G(F)$, we have
 \[
  \deg(\pi)=C\cdot \frac{\dim \rho_\pi}{\#S_{\phi_\pi}}\lvert\gamma(0,\Ad\circ\phi_\pi,\psi)\rvert.
 \]
 Here
 \begin{itemize}
  \item $C\in \R_{>0}$ is a constant which depends only on the Haar measure on $G(F)$ and the non-trivial additive character $\psi$
	(we may take $C=1$ by constructing a Haar measure carefully from $\psi$),
  \item $\Ad\circ \phi_\pi$ is the composite of 
	\[
	 W_F\times\SL_2(\C)\xrightarrow{\phi_\pi}\widehat{G}(\C)\xrightarrow{\Ad}\GL(\Lie \widehat{G}(\C)),
	\]
  \item and 
	\[
	 \gamma(s,\Ad\circ \phi_\pi,\psi)=\varepsilon(s,\Ad\circ \phi_\pi,\psi)\frac{L(1-s,\Ad\circ \phi_\pi)}{L(s,\Ad\circ \phi_\pi)}
	\]
	denotes the local $\gamma$-factor.
 \end{itemize}
\end{conj}

\begin{rem}\label{rem:S-gp-abelian}
 In the case $G=\Sp_{2n}$, $S_{\phi_\pi}$ is known to be an elementary $2$-group, hence $\dim \rho_\pi=1$.
\end{rem}

The formal degree conjecture for a general connected reductive group is formulated similarly,
but one needs slight modification if the group is not simply connected. See \cite{MR2350057} for detail.
The formal degree conjecture has been proved for general linear groups \cite[Theorem 3.1]{MR2350057},
odd special orthogonal groups \cite{MR3649356}, and unitary groups \cite{2018arXiv181200047B},
but it seems still open for many other groups,
such as symplectic groups and even special orthogonal groups.

In this paper, we focus on a very special class of discrete series representations, called simple supercuspidal representations. 
They are constructed by using the compact induction from a compact open subgroup of $G(F)$.
Here we choose a symplectic form given by the skew-symmetric $2n\times 2n$ matrix
\[
 \begin{pmatrix}
  & & & & 1\\ & & & -1 & \\ & & \reflectbox{$\ddots$} & & \\ & 1 & & &\\ -1 & & & &
 \end{pmatrix}
\]
to define $G=\Sp_{2n}$. We define a sequence of compact open subgroups
$G(F)\supset I\rhd I^+\rhd I^{++}$ as follows:
\begin{align*}
 I&=\begin{pmatrix}
    \mathcal{O}_F& &\mathcal{O}_F\\ & \ddots & \\ \mathfrak{p}_F & & \mathcal{O}_F
   \end{pmatrix},\qquad
 I^+=\begin{pmatrix}
    1+\mathfrak{p}_F& &\mathcal{O}_F\\ & \ddots & \\ \mathfrak{p}_F & & 1+\mathfrak{p}_F
   \end{pmatrix},\\
 I^{++}&=\begin{pmatrix}
	 1+\mathfrak{p}_F& \mathfrak{p}_F &  &\mathcal{O}_F\\ & \ddots &\ddots & \\ 
	  & \mathfrak{p}_F & \ddots & \mathfrak{p}_F\\
	  \mathfrak{p}_F^2 & & & 1+\mathfrak{p}_F
	\end{pmatrix}.
\end{align*}
Here $\mathcal{O}_F$ denotes the ring of integers of $F$, and $\mathfrak{p}_F$ the maximal ideal of $\mathcal{O}_F$.
If we fix a uniformizer $\varpi$ of $\mathcal{O}_F$, then we have an isomorphism
\[
 I^+/I^{++}\xrightarrow{\cong} k^{n+1}; (a_{ij})\mapsto (a_{12}\bmod\mathfrak{p}_F,\ldots,a_{n,n+1}\bmod\mathfrak{p}_F,\varpi^{-1}a_{2n,1}\bmod\mathfrak{p}_F).
\]
A character of $I^+/I^{++}\cong k^{n+1}$ is said to be affine generic if it is non-trivial on each factor of $k^{n+1}$.
Let $\chi$ be a character of $\pm I^+$ such that $\chi\vert_{I^{++}}$ is trivial and 
$\chi\vert_{I^+}$ induces an affine generic character of $I^+/I^{++}$.
Then, the compact induction $\cInd_{\pm I^+}^{G(F)} \chi$ is known to be irreducible supercuspidal.
Representations obtained in this way are called simple supercuspidal representations.

The parameter $(\phi_\pi,\rho_\pi)$ attached to a simple supercuspidal representation $\pi$ is investigated by
Oi in detail.

\begin{thm}[{{\cite[Corollary 5.13, Theorem 7.17]{Oi-ssc-classical}}}]\label{thm:ssc-LLC}
 Assume $p\neq 2$. Let $\iota$ denote the embedding $\widehat{G}(\C)=\SO_{2n+1}(\C)\hookrightarrow \GL_{2n+1}(\C)$.
 For a simple supercuspidal representation $\pi$ of $G(F)$,
 we have the following:
 \begin{itemize}
  \item $\iota\circ \phi_\pi=\tau\oplus \omega$, where $\tau$ is an irreducible $2n$-dimensional irreducible representation
	of $W_F$ with Swan conductor $1$ and $\omega$ is a quadratic character of $W_F$.
	Furthermore, $\tau$ is orthogonal, that is, there exists a $W_F$-invariant non-degenerate symmetric bilinear form
	$\tau\times\tau\to \C$.
  \item $\#S_{\phi_\pi}=2$.
 \end{itemize}
\end{thm}

Strictly speaking, \cite[Theorem 7.17]{Oi-ssc-classical} claims that $\iota\circ \phi_\pi=\tau\oplus \det\circ\tau$, 
where $\tau$ is the Langlands parameter of a simple supercuspidal representation of $\GL_n(F)$.
However, it is well-known that such $\tau$ is irreducible and has Swan conductor $1$; see \cite[\S 2]{MR3158004} for example.
Further, since $\tau$ is orthogonal (see \cite[Corollary 5.13]{Oi-ssc-classical}),
the character $\omega=\det\circ\tau$ is quadratic.

By using Theorem \ref{thm:ssc-LLC}, Oi obtained a partial result on the formal degree conjecture.

\begin{thm}[{{\cite[Theorem 9.3]{Oi-ssc-classical}}}]\label{thm:FDC-Oi}
 We write $2n=p^en'$ with $p\nmid n'$.
 Assume $p\neq 2$ and either $p\nmid 2n$ or $n'\mid p-1$.
 Then, Conjecture \ref{conj:HII} holds for simple supercuspidal representations of $G(F)$.
\end{thm}

In the case $p\mid 2n$ and $n'\mid p-1$, Oi used an explicit description of $\tau$ in Theorem \ref{thm:ssc-LLC}
due to Imai and Tsushima \cite{2015arXiv150902960I}, which is extremely complicated.
It involves $4$ field extensions $F^{\mathrm{ur}}\subset E\subset T\subset M\subset N$ of the maximal unramified extension
$F^\mathrm{ur}$ of $F$. The extension $N/E$ is always Galois, but the extension $N/F^{\mathrm{ur}}$ is not necessarily Galois.
The condition $n'\mid p-1$ ensures that $N/F^{\mathrm{ur}}$ is a Galois extension, which makes computations much simpler.

The following is our main theorem for symplectic groups
(for quasi-split even special orthogonal groups, see Remark \ref{rem:Gan-Ichino}).

\begin{thm}\label{thm:main}
 Assume $p\neq 2$. Then, Conjecture \ref{conj:HII} holds for simple supercuspidal representations of $G(F)$.
\end{thm}

In fact, Theorem \ref{thm:main} can be deduced from Theorem \ref{thm:exterior-square} exactly in the same way
as Oi did in \cite[\S 9.3]{Oi-ssc-classical}. We include some of his arguments for reader's convenience. 
In the following, we assume that $p\neq 2$ and let $\pi$ be a simple supercuspidal representation of $G(F)$.
Let $\tau$ and $\omega$ be as in Theorem \ref{thm:ssc-LLC}.

\begin{lem}\label{lem:Artin-L}
 \begin{enumerate}
  \item We have $L(s,\Ad\circ\phi_\pi)=1$.
  \item We have $\Ar(\Ad\circ\phi_\pi)=2n^2+2n$, where $\Ar$ denotes the Artin conductor.
 \end{enumerate}
\end{lem}

\begin{prf}
 First of all, note that $\Ad\circ\phi_\pi=\wedge^2(\tau\oplus \omega)=\wedge^2\tau\oplus \tau\otimes\omega$.
 \begin{enumerate}
  \item It suffices to show that $(\Ad\circ\phi_\pi)^{I_F}=0$. Since $p\neq 2$, the quadratic character $\omega$ is tamely ramified.
	Therefore $(\tau\otimes\omega)^{P_F}=\tau^{P_F}\otimes\omega=0$ by Lemma \ref{lem:totally-wild}.
	Hence it suffices to prove that $(\wedge^2\tau)^{I_F}=0$.

	By Lemma \ref{lem:totally-wild}, $\tau\vert_{I_F}$ is irreducible, hence $\dim (\tau\otimes \tau^\vee)^{I_F}=1$.
	Since $\tau$ is orthogonal by Theorem \ref{thm:ssc-LLC}, we have
	\[
	1=\dim (\tau\otimes \tau^\vee)^{I_F}=\dim (\tau\otimes \tau)^{I_F}=\dim (\Sym^2 \tau)^{I_F}\oplus \dim (\wedge^2 \tau)^{I_F}
	\]
	and $(\Sym^2 \tau)^{I_F}=(\Sym^2 \tau)_{I_F}\neq 0$. Therefore we obtain $(\wedge^2 \tau)^{I_F}=0$, as desired.
  \item Since $\omega$ is tame, we have
	\[
	\Sw(\Ad\circ\phi_\pi)=\Sw(\wedge^2\tau)+\Sw(\tau\otimes\omega)=\Sw(\wedge^2\tau)+\Sw(\tau)=n
	\]
	by Theorem \ref{thm:exterior-square}. On the other hand, by the proof of (i), we have 
	\[
	\dim (\Ad\circ\phi_\pi)/(\Ad\circ\phi_\pi)^{I_F}=\dim (\Ad\circ\phi_\pi)=\frac{2n(2n+1)}{2}=2n^2+n.
	\]
	Therefore we conclude that
	\[
	 \Ar(\Ad\circ\phi_\pi)=\Sw(\Ad\circ\phi_\pi)+\dim (\Ad\circ\phi_\pi)/(\Ad\circ\phi_\pi)^{I_F}=2n^2+2n.
	\]
 \end{enumerate}
\end{prf}

\begin{prf}[of Theorem \ref{thm:main}]
 Let $\St$ denote the Steinberg representation of $G(F)$. 
 We may choose $\psi$ so that the following equalities hold:
 \begin{align*}
  \frac{\deg(\pi)}{\deg(\St)}&=\frac{q^{n^2+n}}{2\gamma(0,\Ad\circ\phi_{\St},\psi)},\\
  \lvert \varepsilon(0,\Ad\circ\phi_\pi,\psi)\rvert&=q^{\frac{1}{2}\Ar(\Ad\circ\phi_\pi)}.
 \end{align*}
 See \cite[(72)]{MR2730575} for the first equality, and \cite[(10) and Proposition 2.3]{MR2730575} for the second.
 Together with Lemma \ref{lem:Artin-L}, we obtain $\lvert \gamma(0,\Ad\circ\phi_\pi,\psi)\rvert=q^{n^2+n}$ and
 \[
  \deg(\pi)=\biggl\lvert \frac{\deg(\St)\gamma(0,\Ad\circ\phi_\pi,\psi)}{2\gamma(0,\Ad\circ\phi_{\St},\psi)}\biggr\rvert.
 \]

 On the other hand, by \cite[\S 3.3]{MR2350057}, the formal degree conjecture for $\St$ is known:
 \[
  \deg(\St)=C\lvert\gamma(0,\Ad\circ\phi_{\St},\psi)\rvert.
 \]
 Hence we have 
 \[
  \deg(\pi)=C\cdot \frac{1}{2}\lvert\gamma(0,\Ad\circ\phi_\pi,\psi)\rvert=C\cdot \frac{\dim\rho_\pi}{\#S_{\phi_\pi}}\lvert\gamma(0,\Ad\circ\phi_\pi,\psi)\rvert,
 \]
 as desired (recall that $\dim \rho_\pi=1$ by Remark \ref{rem:S-gp-abelian} and $\#S_{\phi_\pi}=2$ by Theorem \ref{thm:ssc-LLC}).
\end{prf}

\begin{rem}\label{rem:Gan-Ichino}
 As remarked in \cite[\S 9]{Oi-ssc-classical}, by using the results in \cite{MR3166215}, we can deduce from Theorem \ref{thm:main}
 the formal degree conjecture for simple supercuspidal representations of quasi-split even special orthogonal groups,
 under the assumption $p\neq 2$.
\end{rem}

\def\cftil#1{\ifmmode\setbox7\hbox{$\accent"5E#1$}\else
  \setbox7\hbox{\accent"5E#1}\penalty 10000\relax\fi\raise 1\ht7
  \hbox{\lower1.15ex\hbox to 1\wd7{\hss\accent"7E\hss}}\penalty 10000
  \hskip-1\wd7\penalty 10000\box7}
  \def\cftil#1{\ifmmode\setbox7\hbox{$\accent"5E#1$}\else
  \setbox7\hbox{\accent"5E#1}\penalty 10000\relax\fi\raise 1\ht7
  \hbox{\lower1.15ex\hbox to 1\wd7{\hss\accent"7E\hss}}\penalty 10000
  \hskip-1\wd7\penalty 10000\box7}
  \def\cftil#1{\ifmmode\setbox7\hbox{$\accent"5E#1$}\else
  \setbox7\hbox{\accent"5E#1}\penalty 10000\relax\fi\raise 1\ht7
  \hbox{\lower1.15ex\hbox to 1\wd7{\hss\accent"7E\hss}}\penalty 10000
  \hskip-1\wd7\penalty 10000\box7}
  \def\cftil#1{\ifmmode\setbox7\hbox{$\accent"5E#1$}\else
  \setbox7\hbox{\accent"5E#1}\penalty 10000\relax\fi\raise 1\ht7
  \hbox{\lower1.15ex\hbox to 1\wd7{\hss\accent"7E\hss}}\penalty 10000
  \hskip-1\wd7\penalty 10000\box7} \def\cprime{$'$} \def\cprime{$'$}
  \newcommand{\dummy}[1]{}
\providecommand{\bysame}{\leavevmode\hbox to3em{\hrulefill}\thinspace}
\providecommand{\MR}{\relax\ifhmode\unskip\space\fi MR }
\providecommand{\MRhref}[2]{%
  \href{http://www.ams.org/mathscinet-getitem?mr=#1}{#2}
}
\providecommand{\href}[2]{#2}

\end{document}